\documentclass[a4paper,12pt,fleqn]{article}
\usepackage[latin1]{inputenc}  \usepackage[T1]{fontenc}
\usepackage{lscape}
\usepackage{tabularx}     
\usepackage{pstricks}    
\usepackage{rotating}
\usepackage{amsfonts,amsmath,amsthm,amssymb,dsfont}
\usepackage{marginnote}
\usepackage{rotate}
\usepackage{pgfplots}
\usepackage{tikz}
\usepackage{tikz-cd}
\usepackage[arrow,matrix,curve]{xy} 	
\usepackage{xcolor}
\title{Deducing Three Gap Theorem From Rauzy-Veech Induction}
\author{Christian Weiß}
\date{\today}

\newtheorem{thm}{Theorem}

\newcommand{\RR}{{\mathbb{R}}}
\newcommand{\ZZ}{{\mathbb{Z}}}
\newcommand{\NN}{{\mathbb{N}}}

\makeindex

\begin{document} 

\maketitle

\begin{abstract}% englische Fassung
	The Three Gap Theorem states that there are at most three distinct lengths of gaps if one places $n$ points on a circle, at angles of $z, 2z, 3z, \ldots nz$ from the starting point. The theorem was first proven in 1958 by S\'os and many proofs have been found since then. In this note we show how the Three Gap Theorem can easily be deduced by using Rauzy-Veech induction.
\end{abstract}

%\tableofcontents

Kronecker sequences are important examples of uniformly distributed sequences. Given $z \in \RR$, a Kronecker sequence is of the form $(z_n)_{n \geq 0} = (\left\{ nz \right\})_{n \geq 0}$ where $\left\{ z \right\} := z - \lfloor z \rfloor$ denotes the fractional part of $z$.  Sometimes their level of uniformity is even as great as possible since infinitely many Kronecker sequences belong to the classical examples of low-discrepancy sequences, see e.g. \cite{DT97}, Corollary~1.65. Another important property describing the uniformity of Kronecker sequences is the Three Gap Theorem, going back to a famous conjecture of Steinhaus which was first proved by S\'os in \cite{Sos58}. Since the Three Gap Theorem is closely linked to continued fraction expansion we shortly introduce some notation first: every irrational number $z$ has a uniquely determined infinite continued fraction expansion
$$z = a_0 + 1/(a_1+1/(a_2+\ldots)) =: [a_0;a_1;a_2;\ldots],$$
where the $a_i$ are integers with $a_0 = \lfloor z \rfloor$ and $a_i \geq 1$ for all $i \geq 1$. The sequence of \textbf{convergents} $(r_i)_{i \in \NN}$ of $z$ is defined by
$$r_i = [a_0;a_1;\ldots;a_i].$$
The convergents $r_i = p_i/q_i$ with $\gcd(p_{i},q_{i}) = 1$ can also be calculated directly by the recurrence relation
\begin{align} \label{eqs1}
\begin{split}
& p_{-2} = 0, \qquad p_{-1} = 1, \qquad p_{i} = a_ip_{i-1} + p_{i-2}, \quad i \geq 0\\
& q_{-2} = 1, \qquad q_{-1} = 0, \qquad q_{i} = a_iq_{i-1} + q_{i-2}, \quad i \geq 0.
\end{split}
\end{align}
Now the Three Gap Theorem states the following.
\begin{thm}[Three Gap Theorem] \label{thm:3gap} Let $z \in (0,1)$ be irrational with continued fraction expansion $z =[a_0;a_1;a_2;\ldots]$ and convergents $r_n = p_n/q_n$. Furthermore let $N \in \NN$ with $N \geq 2$ have Ostrowski representation\footnote{Actually, this is a slightly amended version of the Ostrowski representation, since usually it is assumed that $0 \leq b_0 < a_1$ and $b_{j-1} = 0$ if $b_j = a_{j+1}$ but not $q_m + 1\leq N < q_{m+1} + q_m$.}
	$$N = \sum_{j=0}^m b_j q_j$$
with integer coefficients $b_0,b_1,\ldots,b_m$ satisfying $0 \leq b_j \leq a_{j_+1}$ and minimal $m$ with $q_m + 1\leq N < q_{m+1} + q_m$. If $ z < \tfrac{1}{2}$, let $K_{2l-1} = \left\{ q_{2l-1} z \right\}$ and $K_{2l} = 1 - \left\{ q_{2l} z \right\}$ for $l \in \NN_0$, and if $z > \tfrac{1}{2}$, let $K_{2l-1} = 1- \left\{ q_{2l-1} z \right\}$ and $K_{2l} = \left\{ q_{2l} z \right\}$. Then the finite sequence $(\left\{ n z \right\})_{n=1,\ldots,N-1}$ has at most three different lengths of gaps, namely %	L_1 & = \max(K_{m},K_{m-1}) - b_m \min(K_{m},K_{m-1}), \\L_2 & = \min(K_{m},K_{m-1}),\\
\begin{align*}
	L_1 & = K_{m-1} - b_m K_m, \\
	L_2 & = K_m\\
	L_3 & = L_1 + L_2.
\end{align*}
The number of gaps of lengths $L_1, L_2, L_3$ are
\begin{align*}
	N_1 & = N - b_m q_m - q_{m-1}, \\
	N_2 & = N-q_{m},\\
	N_3 & = q_m - (N - b_m q_m - q_{m-1}).
\end{align*}
 \end{thm}
Later on, several further proofs of the claim have been found, see e.g. \cite{Lia79}, \cite{AB98}, \cite{MS17} and most recently \cite{Tah17}. In this note we add another proof to the list of proofs of Theorem~\ref{thm:3gap}: a Kronecker sequence corresponds to the rotation of the unit circle by the angle $2 \pi z$ if we identify $[0,1)$ with $\RR / \ZZ$. The map $x \mapsto x + z \mod 1$ is also the simplest example of an interval exchange transformation. Indeed, it may be considered an exchange of two subintervals of $[0,1)$, specifically $A:=[0,1-\left\{z\right\})$ and $B:=[1-\left\{z\right\},1).$ Here we show that the Three Gap Theorem can easily be deduced from the viewpoint of interval exchange transformations by using Rauzy-Veech induction.

\paragraph{Interval Exchange Transformations.} We only give a brief summary here and refer the reader to \cite{Via06}, which is also our main source, for more details. Let $I \subset \RR$ be an interval of the form $[0,\lambda^*)$ and let $\left\{ I_\alpha | \alpha \in \mathcal{A} \right\}$ be a finite partition of $I$ into subintervals indexed by some finite alphabet $\mathcal{A}$. An \textbf{interval exchange transformation} is a map $f: I \to I$ which is a translation on each subinterval $I_\alpha$. It is determined by its combinatorial data and its length data. The \textbf{combinatorial data} consists of two bijections $\pi_0, \pi_1: \mathcal{A} \to \left\{1,\ldots,n \right\}$, where $n$ is the number of elements of $\mathcal{A}$ and the \textbf{length data} are numbers $(\lambda_\alpha)_{\alpha \in \mathcal{A}}=:\lambda$ with $\lambda_\alpha > 0$ and $\lambda^* = \sum_{\alpha \in \mathcal{A}} \lambda_\alpha$. The number $\lambda_\alpha$ is the length of the subinterval $I_\alpha$ and the pair $\pi = (\pi_0,\pi_1)$ describes the ordering of the subintervals before and after the map $f$ is iterated. For $\mathcal{A} = \left\{A, B \right\}, \pi_0(A) = \pi_1(B) = 1$ and $\pi_1(A) = \pi_0(B) = 2$ the interval exchange transformation becomes the rotation of $\RR / \lambda^* \ZZ$ by $\lambda_B$.%In analytical terms, consider for $\epsilon \in \left\{ 0, 1 \right\}$ the maps $j_\epsilon$ which are given on $I_\alpha$ by
  %\begin{align} \label{eq5}
  %j_\epsilon(x) = x + \sum_{\pi_\epsilon(\beta) < \pi_\epsilon(\alpha)} \lambda_\beta.
  %\end{align}
  %The interval exchange transformation corresponding to the combinatorial and length data equals $f = j_1 \circ j_0^{-1}$. 
  
 \paragraph{Rauzy-Veech Induction.} For an interval exchange transformation $f$ given by the data $(\pi,\lambda)$, we define the \textbf{type $\epsilon$} by 
 $$\lambda_{\pi_\epsilon^{-1}(n)} >  \lambda_{\pi_{1-\epsilon}^{-1}(n)},$$
 if the lengths of these two intervals do not coincide.  Following the usual terminology, we say that $I_{\pi_{\epsilon}^{-1}(n)}$ is the \textbf{winner} and $I_{\pi_{1-\epsilon}^{-1}(n)}$ is the \textbf{loser}. Let $J$ be the subinterval of $I$ obtained by removing the loser, i.e.
 $$J = \begin{cases} I \setminus f(I_{\pi^{-1}_1(n)}) & \textrm{if $f$ has type } 0 \\ I \setminus I_{\pi^{-1}_0(n)} & \textrm{if $f$ has type } 1. \end{cases}$$
 The \textbf{Rauzy-Veech induction} $R(f)$ of $f$ is its first return map to the interval $J$. It is again an interval exchange transformation consisting of $n$ subintervals. The corresponding data $(\pi',\lambda')$ can be easily calculated, see e.g. \cite{Via06}. In this note, we restrict our attention to $\pi_0(A) = \pi_1(B) = 1$ and $\pi_1(A) = \pi_0(B) = 2$ because it is the only case of interest for the proof of Theorem~\ref{thm:3gap}. The first $a-1$ steps of Rauzy-Veech induction for two intervals are depicted in Figure~1, where $A^{k}$ respectively $B^{k}$ denotes the specific interval appearing after step $k$.\\[12pt]
 As long as the two rightmost intervals $\lambda_{\pi_\epsilon^{-1}(n)}, \lambda_{\pi_{1-\epsilon}^{-1}(n)}$ have different lengths the construction can be iterated yielding a sequence $R^{(j)}(f) = (\pi^{(j)},\lambda^{(j)})$ of interval exchange transformations. The type is well-defined infinitely times if and only if the so-called \textbf{Keane-condition} which postulates that the orbits of the singularities of $f^{-1}$ by $f$ are infinite and distinct is satisfied, see e.g \cite{Yoc06}.
 \begin{center}
 	\begin{tikzpicture}[scale=0.75]
 	\def\x{10};
 	\def\klein{0.25};
 	\def\gross{1};
 	
 	%links 	
 	\draw[thick] (0,0)--(8,0);     	
 	\draw[thick] (0,0)--(0,0.25);     	
 	\draw (0.6,0)--(0.6,0) node[above] {\scriptsize{$A$}};
 	\draw[thick] (1.2,0)--(1.2,0.25);
 	\draw (4.6,0)--(4.6,0) node[above] {\scriptsize{$B$}};    	
 	\draw[thick] (8,0)--(8,0.25);

 	\draw[thick] (0,0-\klein)--(8,0-\klein);     	
 	\draw[thick] (0,0-\klein)--(0,-0.25-\klein);     	
 	\draw (3.4,0-\klein)--(3.4,0-\klein) node[below] {\scriptsize{$B$}};
 	\draw[thick] (6.8,0-\klein)--(6.8,-0.25-\klein);
 	\draw (7.4,0-\klein)--(7.4,0-\klein) node[below] {\scriptsize{$f(A)$}};    	
 	\draw[thick] (8,0-\klein)--(8,-0.25-\klein);
 	
 	\draw[thick,->] (4,-\gross)--(4,-\gross - 0.5);

 	\draw[thick] (0,-2*\gross-\klein)--(6.8,-2*\gross-\klein);
 	\draw[thick,dashed,lightgray] (6.8,-2*\gross-\klein)--(8,-2*\gross-\klein);  
 	\draw[thick] (0,-2*\gross-\klein)--(0,-2*\gross);     	
 	\draw (0.6,-2*\gross-\klein)--(0.6,-2*\gross-\klein) node[above] {\tiny{$A$}};    
 	\draw[thick] (1.2,-2*\gross-\klein)--(1.2,-2*\gross);
 	\draw (4.0,-2*\gross-\klein)--(4.0,-2*\gross-\klein) node[above] {\tiny{$B^{1}$}};
 	\draw[thick] (6.8,-2*\gross-\klein)--(6.8,-2*\gross);
 	\draw (7.4,-2*\gross-\klein)--(7.4,-2*\gross-\klein) node[above,lightgray] {\tiny{$f(A)$}};    	
 	\draw[thick,lightgray] (8,-2*\gross-\klein)--(8,-2*\gross);
 	
 	\draw[thick] (0,-2*\gross-2*\klein)--(6.8,-2*\gross-2*\klein);
 	\draw[thick,dashed,lightgray] (6.8,-2*\gross-2*\klein)--(8,-2*\gross-2*\klein);  
 	\draw[thick] (0,-2*\gross-2*\klein)--(0,-2*\gross-3*\klein);     	
 	\draw (2.8,-2*\gross-2*\klein)--(2.8,-2*\gross-2*\klein) node[below] {\tiny{$B^{1}$}};    
 	\draw[thick] (5.6,-2*\gross-2*\klein)--(5.6,-2*\gross-3*\klein);
 	\draw (6.2,-2*\gross-2*\klein)--(6.2,-2*\gross-2*\klein) node[below] {\tiny{$f^{2}(A)$}};
 	\draw[thick] (6.8,-2*\gross-2*\klein)--(6.8,-2*\gross-3*\klein);
 	\draw (7.4,-2*\gross-2*\klein)--(7.4,-2*\gross-2*\klein) node[below,lightgray] {\tiny{$f(A)$}};    	
 	\draw[thick,lightgray] (8,-2*\gross-2*\klein)--(8,-2*\gross-3*\klein);
 	
 	\draw[thick,->] (4,-3*\gross)--(4,-3*\gross-2*\klein);
 	
 	\draw[thick,dotted] (3.5,-3*\gross - 3*\klein)--(4.5,-3*\gross - 3*\klein);
 	
 	\draw[thick,->] (4,-4*\gross)--(4,-4*\gross-2*\klein);
 	
 	\draw[thick] (0,-5*\gross-\klein)--(2.0,-5*\gross-\klein);
 	\draw[thick,dashed,lightgray] (2.0,-5*\gross-\klein)--(8,-5*\gross-\klein);  
 	\draw[thick] (0,-5*\gross-\klein)--(0,-5*\gross);     	
 	\draw (0.6,-5*\gross-\klein)--(0.6,-5*\gross-\klein) node[above] {\tiny{$A$}};    
 	\draw[thick] (1.2,-5*\gross-\klein)--(1.2,-5*\gross);
 	\draw (1.6,-5*\gross-\klein)--(1.6,-5*\gross-\klein) node[above] {\tiny{$B^{a-1}$}};
 	\draw[thick] (2.0,-5*\gross-\klein)--(2.0,-5*\gross); 	
 	\draw (3.8,-5*\gross-\klein)--(3.8,-5*\gross-\klein) node[above,lightgray] {\ldots};  
 	\draw[thick,lightgray] (5.6,-5*\gross-\klein)--(5.6,-5*\gross);
 	\draw (6.2,-5*\gross-\klein)--(6.2,-5*\gross-\klein) node[above,lightgray] {\tiny{$f^{2}(A)$}};  	  	
 	\draw[thick,lightgray] (6.8,-5*\gross-\klein)--(6.8,-5*\gross);
 	\draw (7.4,-5*\gross-\klein)--(7.4,-5*\gross-\klein) node[above,lightgray] {\tiny{$f(A)$}};    	
 	\draw[thick,lightgray] (8,-5*\gross-\klein)--(8,-5*\gross);
 	
 	\draw[thick] (0,-5*\gross-2*\klein)--(2.0,-5*\gross-2*\klein);
 	\draw[thick,dashed,lightgray] (2.0,-5*\gross-2*\klein)--(8,-5*\gross-2*\klein);  
 	\draw[thick] (0,-5*\gross-2*\klein)--(0,-5*\gross-3*\klein);     	
 	\draw (0.4,-5*\gross-2*\klein)--(0.4,-5*\gross-2*\klein) node[below] {\tiny{$B^{a-1}$}};    
 	\draw[thick] (0.8,-5*\gross-2*\klein)--(0.8,-5*\gross-3*\klein);
 	\draw (1.4,-5*\gross-2*\klein)--(1.4,-5*\gross-2*\klein) node[below] {\tiny{$f^{a}(A)$}};
 	\draw[thick] (2.0,-5*\gross-2*\klein)--(2.0,-5*\gross-3*\klein);
 	\draw(3.8,-5*\gross-2*\klein)--(3.8,-5*\gross-2*\klein) node[below,lightgray] {\ldots};  
 	\draw[thick,lightgray] (5.6,-5*\gross-2*\klein)--(5.6,-5*\gross-3*\klein);
 	\draw (6.2,-5*\gross-2*\klein)--(6.2,-5*\gross-2*\klein) node[below,lightgray] {\tiny{$f^{2}(A)$}};  	  	
 	\draw[thick,lightgray] (6.8,-5*\gross-2*\klein)--(6.8,-5*\gross-3*\klein);
 	\draw (7.4,-5*\gross-2*\klein)--(7.4,-5*\gross-2*\klein) node[below,lightgray] {\tiny{$f(A)$}};    	
 	\draw[thick,lightgray] (8,-5*\gross-2*\klein)--(8,-5*\gross-3*\klein);

 	%rechts
 	
 	\draw[thick] (0+\x,0)--(8+\x,0);     	
 	\draw[thick] (0+\x,0)--(0+\x,0.25);     	
 	\draw (3.4+\x,0)--(3.4+\x,0) node[above] {\scriptsize{$A$}};
 	\draw[thick] (6.8+\x,0)--(6.8+\x,0.25);
 	\draw (7.4+\x,0)--(7.4+\x,0) node[above] {\scriptsize{$B$}};    	
 	\draw[thick] (8+\x,0)--(8+\x,0.25);

 	\draw[thick] (0+\x,0-\klein)--(8+\x,0-\klein);     	
 	\draw[thick] (0+\x,0-\klein)--(0+\x,-0.25-\klein);     	
 	\draw (0.6+\x,0-\klein)--(0.6+\x,0-\klein) node[below] {\scriptsize{$f(B)$}};
 	\draw[thick] (1.2+\x,0-\klein)--(1.2+\x,-0.25-\klein);
 	\draw (4.6+\x,0-\klein)--(4.6+\x,0-\klein) node[below] {\scriptsize{$A$}};    	
 	\draw[thick] (8+\x,0-\klein)--(8+\x,-0.25-\klein);
 	
 	\draw[thick,->] (4+\x,-\gross)--(4+\x,-\gross - 0.5);

 	\draw[thick] (0+\x,-2*\gross-\klein)--(6.8+\x,-2*\gross-\klein);
 	\draw[thick,dashed,lightgray] (6.8+\x,-2*\gross-\klein)--(8+\x,-2*\gross-\klein);  
 	\draw[thick] (0+\x,-2*\gross-\klein)--(0+\x,-2*\gross);     	
 	\draw (2.8+\x,-2*\gross-\klein)--(2.8+\x,-2*\gross-\klein) node[above] {\tiny{$A^{1}$}};    
 	\draw[thick] (5.6+\x,-2*\gross-\klein)--(5.6+\x,-2*\gross);
 	\draw (6.2+\x,-2*\gross-\klein)--(6.2+\x,-2*\gross-\klein) node[above] {\tiny{$f^{-1}(B)$}};
 	\draw[thick] (6.8+\x,-2*\gross-\klein)--(6.8+\x,-2*\gross);
 	\draw (7.4+\x,-2*\gross-\klein)--(7.4+\x,-2*\gross-\klein) node[above,lightgray] {\tiny{$B$}};    	
 	\draw[thick,lightgray] (8+\x,-2*\gross-\klein)--(8+\x,-2*\gross);
 	
 	\draw[thick] (0+\x,-2*\gross-2*\klein)--(6.8+\x,-2*\gross-2*\klein);
 	\draw[thick,dashed,lightgray] (6.8+\x,-2*\gross-2*\klein)--(8+\x,-2*\gross-2*\klein);  
 	\draw[thick] (0+\x,-2*\gross-2*\klein)--(0+\x,-2*\gross-3*\klein);     	
 	\draw (0.6+\x,-2*\gross-2*\klein)--(0.6+\x,-2*\gross-2*\klein) node[below] {\tiny{$f(B)$}};    
 	\draw[thick] (1.2+\x,-2*\gross-2*\klein)--(1.2+\x,-2*\gross-3*\klein);
 	\draw (4.4+\x,-2*\gross-2*\klein)--(4.4+\x,-2*\gross-2*\klein) node[below] {\tiny{$A^{1}$}};
 	\draw[thick] (6.8+\x,-2*\gross-2*\klein)--(6.8+\x,-2*\gross-3*\klein);
 	\draw (7.4+\x,-2*\gross-2*\klein)--(7.4+\x,-2*\gross-2*\klein) node[below,lightgray] {\tiny{$B$}};    	
 	\draw[thick,lightgray] (8+\x,-2*\gross-2*\klein)--(8+\x,-2*\gross-3*\klein);
 	
 	\draw[thick,->] (4+\x,-3*\gross)--(4+\x,-3*\gross-2*\klein);
 	
 	\draw[thick,dotted] (3.5+\x,-3*\gross - 3*\klein)--(4.5+\x,-3*\gross - 3*\klein);
 	
 	\draw[thick,->] (4+\x,-4*\gross)--(4+\x,-4*\gross-2*\klein);
 	
 	\draw[thick] (0+\x,-5*\gross-\klein)--(2.0+\x,-5*\gross-\klein);
 	\draw[thick,dashed,lightgray] (2.0+\x,-5*\gross-\klein)--(8+\x,-5*\gross-\klein);  
 	\draw[thick] (0+\x,-5*\gross-\klein)--(0+\x,-5*\gross);     	
 	\draw (0.4+\x,-5*\gross-\klein)--(0.4+\x,-5*\gross-\klein) node[above] {\tiny{$A^{a-1}$}};    
 	\draw[thick] (0.8+\x,-5*\gross-\klein)--(0.8+\x,-5*\gross);
 	\draw (1.4+\x,-5*\gross-\klein)--(1.4+\x,-5*\gross-\klein) node[above] {\tiny{$f^{-a}(B)$}};
 	\draw[thick] (2.0+\x,-5*\gross-\klein)--(2.0+\x,-5*\gross); 	
 	\draw (3.8+\x,-5*\gross-\klein)--(3.8+\x,-5*\gross-\klein) node[above,lightgray] {\ldots};  
 	\draw[thick,lightgray] (5.6+\x,-5*\gross-\klein)--(5.6+\x,-5*\gross);
 	\draw (6.2+\x,-5*\gross-\klein)--(6.2+\x,-5*\gross-\klein) node[above,lightgray] {\tiny{$f^{-1}(B)$}};  	  	
 	\draw[thick,lightgray] (6.8+\x,-5*\gross-\klein)--(6.8+\x,-5*\gross);
 	\draw (7.4+\x,-5*\gross-\klein)--(7.4+\x,-5*\gross-\klein) node[above,lightgray] {\tiny{$B$}};    	
 	\draw[thick,lightgray] (8+\x,-5*\gross-\klein)--(8+\x,-5*\gross);
 	
 	\draw[thick] (0+\x,-5*\gross-2*\klein)--(2.0+\x,-5*\gross-2*\klein);
 	\draw[thick,dashed,lightgray] (2.0+\x,-5*\gross-2*\klein)--(8+\x,-5*\gross-2*\klein);  
 	\draw[thick] (0+\x,-5*\gross-2*\klein)--(0+\x,-5*\gross-3*\klein);     	
 	\draw (0.6+\x,-5*\gross-2*\klein)--(0.6+\x,-5*\gross-2*\klein) node[below] {\tiny{$f(B)$}};    
 	\draw[thick] (1.2+\x,-5*\gross-2*\klein)--(1.2+\x,-5*\gross-3*\klein);
 	\draw (1.6+\x,-5*\gross-2*\klein)--(1.6+\x,-5*\gross-2*\klein) node[below] {\tiny{$A^{a-1}$}};
 	\draw[thick] (2.0+\x,-5*\gross-2*\klein)--(2.0+\x,-5*\gross-3*\klein);
 	\draw(3.8+\x,-5*\gross-2*\klein)--(3.8+\x,-5*\gross-2*\klein) node[below,lightgray] {\ldots};  
 	\draw[thick,lightgray] (5.6+\x,-5*\gross-2*\klein)--(5.6+\x,-5*\gross-3*\klein);
 	\draw (6.2+\x,-5*\gross-2*\klein)--(6.2+\x,-5*\gross-2*\klein) node[below,lightgray] {\tiny{$f^{-1}(B)$}};  	  	
 	\draw[thick,lightgray] (6.8+\x,-5*\gross-2*\klein)--(6.8+\x,-5*\gross-3*\klein);
 	\draw (7.4+\x,-5*\gross-2*\klein)--(7.4+\x,-5*\gross-2*\klein) node[below,lightgray] {\tiny{$B$}};    	
 	\draw[thick,lightgray] (8+\x,-5*\gross-2*\klein)--(8+\x,-5*\gross-3*\klein);
 	
 	\end{tikzpicture}\\[6pt]
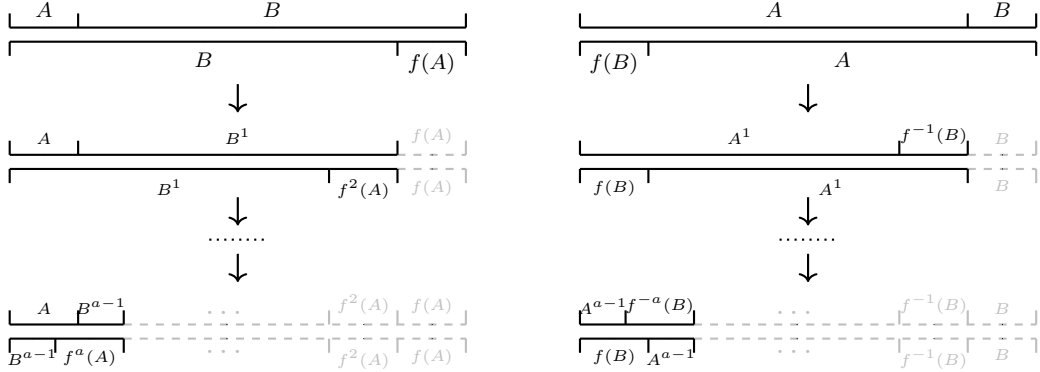
 	Figure 1. First $a$ steps of Rauzy-Veech induction, for type $0$ (left) and type $1$ (right) \label{fig:Rauzy}
 \end{center}
 Besides the Rauzy-Veech induction we need the \textbf{accelerated Rauzy-Veech induction}, also known as \textbf{Zorich-Rauzy-Veech induction}: it means that the Rauzy-Veech induction is applied as many times $a \in \NN$ until the type changes, compare again Figure~1. We denote the resulting map by $\hat{R}$ and the corresponding $a$ by $a_i$. In the case of two intervals, the accelerated Rauzy-Veech induction is equivalent to the continued fraction algorithm and the $a_i$ are equal to the coefficients of the continued fraction expansion, see e.g. \cite{Via06}, Chapter 9.
 
 \paragraph{Proof of the Three Gap Theorem.} Let $N \in \NN$ be arbitrary and let us assume that $1 > z > \tfrac{1}{2}$ since the case $z < \tfrac{1}{2}$ works likewise. At first, we prove by induction that the $m^{th}$ application of \textit{accelerated} Rauzy-Veech yields a partition of $[0,1)$ into $q_m$ long and $q_{m-1}$ short intervals. Using the notation therein, we can see from Figure~1 that the first application of accelerated Rauzy-Veech disjointly partitions
 \begin{align} \label{eq:decomp}
 [0,1) = f(A) \dot{\cup} f^2(A) \dot{\cup} \ldots \dot{\cup} f^{a_1}(A) \dot{\cup} B^{a-1}
 \end{align}
 into $1 = q_0$ short interval with length $\hat{l}_s^{(1)}$ and $a_1  = q_1$ long intervals with length $\hat{l}_l^{(1)}$. In the same manner, the $m$-th application of accelerated Rauzy-Veech partitions the long interval (before $m$-th step) of length $\hat{l}_{l}^{(m-1)}$ into $a_m$ long intervals (after the $m$-th step) of length $\hat{l}_l^{(m)}=\hat{l}_{s}^{(m-1)}$ and one short interval of length $\hat{l}_s^{(m)}$. Inductively it follows from \eqref{eq:decomp}, that applying $f$ to one distinguished long interval $a_m \cdot q_{m-1} + q_{m-2} = q_m$ times (compare \eqref{eqs1}) and to one distinguished small interval $q_{m-1}$ times yields again a disjoint partition of $[0,1)$.\\[12pt]
 Hence, if $N= q_m + q_{m-1}$, then the set of left endpoints of the partition implied by $\hat{R}^m$, which consists of $N = q_m + q_{m-1}$ subintervals, only has two different gap lengths, i.e. $N_1 = 0$, $N_2 = q_{m-1}$ and $N_3 = q_m$, as claimed. Similarly, one step of the usual (not accelerated) Rauzy-Veech algorithm partitions the long interval of length $l_{l}^{(i)}$ into one interval of length $l_{l}^{(i)} - l_{s}^{(i)}$ and one of length $l_s^{(i)}$, compare again Figure~1. Note that $\hat{l}^{(m)} = l^{(\sum_{k=0}^m a_k)}.$ This proves the formulas for $N_1, N_2, N_3$ in the case $N= b_mq_m + q_{m-1}$ since applying the usual Rauzy-Veech algorithm at this stage increases $N$ by $q_m$.\\[12pt]
 Moreover by applying \eqref{eq:decomp} inductively, we see that the set of left endpoints equals the finite Kronecker sequence  $((n-2)z)_{n=1,\ldots,N}$ as a set. This transfers the results on the number of gaps to the Kronecker sequence for $N = b_mq_m + q_{m-1}$. Note that considering  $((n-2)z)_{n=1,\ldots,N}$  instead of  $(nz)_{n=1,\ldots,N-1}$ is no restriction because the variable shift just corresponds to a rotation of all points of the sequence by a fixed angle, namely $-2z$, which does not change the number of gaps or their lengths.\\[12pt] 
 If $i$ Rauzy-Veech steps have been applied and $N$ is increased by $1$, the Kronecker sequence must approach the set of left endpoints implied by the next Rauzy-Veech step $R^{(i+1)}$ as a set. In other words, \textit{one} of the long intervals is subdivided into one subinterval of length $l_s^{(i)}$ and one subinterval of length $l_m := l_l^{(i)} - l_s^{(i)}$ by the additional point of the Kronecker sequence while the other intervals remain undivided. In other words, as $N$ increases by $1$ also $N_1$ (number of medium length intervals) and $N_2$ (number of short intervals) increase by $1$, while $N_3$ (number of long intervals) decreases by $1$ until there is no interval of length $l_l^{(i)}$ left (and thus $N_1$ drops to $0$ again). This completes the proof of the expressions for $N_1,N_2$ and $N_3$ for \textit{all} values of $N$.\\[12pt]
 The only part of the assertion which is missing are the formulas for $L_1, L_2$ and $L_3$. Of course, the total length of all subintervals has to sum up to $1$ and we have already seen $L_3 = L_1  + L_2$ because $l_m = l_l^{(i)} - l_s^{(i)}$. Therefore, it suffices to only calculate one of the lengths. Another time we use induction to show the formula for the shortest length $L_2$. At the beginning there is $1=q_0$ interval of short length $\hat{l}_s^{(0)} = L_2^{(0)} = 1 - \left\{ q_0 z \right\}$ and applying accelerated Rauzy-Veech once yields 
 %$$L_2^{(q_0+q_1)} = \left\{ q_0 z \right\} - a_1 (1 - \left\{ q_0 z\right\}) = 1 - \left\{ q_1 z \right\}.$$
 $$\hat{l}_s^{(1)} = L_2^{(1)} = 1 - a_1 (1 - \left\{ q_0 z\right\}) = \left\{ q_1 z \right\}.$$
 Again by accelerated Rauzy-Veech, the $m$-th appearing short interval length $\hat{l}_s^{(m)} = L_2^{(m)}$ has to satisfy the equation
 \begin{align*}
 	L_2^{(m)} & = L_3^{(m-1)} - a_m \cdot L_2^{(m-1)}\\
 				& = K_{m-2} - a_m K_{m-1}\\
 				& = \begin{cases} (1 - \left\{ q_{m-2}z\right\}) - a_m  \left\{ q_{m-1}z \right\} & \textrm{m is even} \\  \left\{ q_{m-2}z \right\} - a_m (1 - \left\{ q_{m-1}z\right\}) & \textrm{m is odd} \end{cases}\\
 				& = \begin{cases}  1 - \left\{q_mz\right\}& \textrm{m is even}\\ \ \left\{ q_mz \right\} & \textrm{m is odd} \end{cases}
 \end{align*}
 This finishes the proof because the length $L_2$ does not change before one accelerated Rauzy-Veech step is completed.
 %It will turn out to be useful to consider the shifted Kronecker sequence $((n-2)z)_{n=1,\ldots,N}$ which we may without loss of generality do as it just corresponds to a rotation of the sequence $(nz)$.
  
\paragraph{Acknowledgement.} The author thanks two anonymous referees for their very valuable comments.

\textsc{Christian Wei\ss, Hochschule Ruhr West, Duisburger Str. 100, D-45479 M\"ulheim an der Ruhr}\\
\textit{E-mail address:} \texttt{christian.weiss@hs-ruhrwest.de}

\end{document}